\documentclass[11pt,abstracton]{scrartcl}
 \usepackage{amsmath,amsthm,amsfonts,amssymb,bbm,bm}
\usepackage[numbers, square]{natbib}
 \usepackage[pagebackref=true,linktoc=page,colorlinks=true,linkcolor=blue,citecolor=blue]{hyperref}

\usepackage{enumerate}
\usepackage{multirow} %% in tabellen
\usepackage{graphicx}
\usepackage[all]{xy}
\usepackage{tabularx} % displaystyle in tabulars
\usepackage{booktabs} % for toprule/midrule etc
\usepackage{accents} % for accentset for interior
\usepackage[font=small,labelfont=bf]{caption} % small captions

\newcommand{\NN}{\ensuremath{\mathbb{N}}}

\def\ci{\perp\!\!\!\perp}

\theoremstyle{plain}

\theoremstyle{definition}

\theoremstyle{remark}
\newtheorem*{remark}{Remark}
\begin{document}

\title{Extremal independence old and new}

\subtitle{An extended version of the contribution to the RSS discussion\\ on the forthcoming paper ``Graphical models for extremes''\\ by Sebastian~Engelke and Adrien~S.~Hitz}

\author{
  Kirstin Strokorb\footnote{School of Mathematics, Cardiff University, Cardiff CF24 4AG, Email: StrokorbK@cardiff.ac.uk} }

\maketitle 

\begin{abstract}
  On 12~February~2020 the Royal Statistical Society hosted a meeting to discuss the forthcoming paper ``Graphical models for extremes'' by Sebastian Engelke and Adrien Hitz. This short note is a supplement to my discussion contribution. It contains the proofs. It is shown that the traditional notion of extremal indepdendence agrees with the newly introduced notion of extremal independence, which subsequently allows for a meaningful interpretation of disconnected graphs in the context of the discussion paper. The notation and references used in this note are adopted from \cite{eh2020}.
\end{abstract}
\vspace*{1cm}

In the setting of the discussion paper \cite{eh2020} all extremal graphical models have to be \emph{connected} as the exponent measure $\Lambda$ has been chosen to have a Lebesgue-density on the nonnegative orthant $\mathcal{E}=[{0},{\infty})^d\setminus \{\bm{0}\}$. Indeed, as has been pointed out already in \cite{eh2020}, the latter implies extremal dependence between all components of $\bm{X}$ and that the definition of conditional extremal independence $\bm{Y}_A \perp_e \bm{Y}_C | \bm{Y}_B$ cannot be extended to the case where $B=\emptyset$.

  For completeness, let us revisit the broader framework in which the exponent measure $\Lambda$ is allowed to place mass on all faces
  \begin{align*}
    \mathcal{E}^I=\{\bm{x} \in  \mathcal{E} : \bm{x}_I > \bm{0}, \bm{x}_{\setminus I} = \bm{0}\} \qquad
\text{of} \qquad \mathcal{E} = \bigcup_{\emptyset \neq I \subset \{1,\dots,d\}} \mathcal{E}^I
  \end{align*}
and let $\{1,\dots,d\}=A\cup C$ with $A \cap C = \emptyset$. In this setting, the analog of conditional extremal independence (17) for $B=\emptyset$ (the unconditional case) is
\begin{align*}
  \bm{Y}^k_A \ci \bm{Y}^k_C \quad \text{for all  $k=1,\dots,d$}. 
\end{align*}
It defines a \textit{new notion of extremal independence} that we may abbreviate by $\bm{Y}_A \perp_e \bm{Y}_C$.

Naturally, the question arises how it is related to the \textit{traditional notion of extremal independence} (which can be expressed in terms of the support set of the exponent measure $\Lambda$, see below). The answer is as follows and its proof does not require any considerations on densities.

\paragraph{Theorem.} \textit{The new notion of extremal independence and the traditional notion of extremal indepdendence are equivalent.}\\

\newpage

\vspace*{1cm}
To be clear about the traditional notion of extremal independence, it is convenient to abbreviate $\mathcal{E}_I=[{0},{\infty})^I \setminus \{\bm{0}\}$ and for $\emptyset \neq J \subset I \subset \{1,\dots,d\}$
  \begin{align*}\mathcal{E}^J_I=\{\bm{x} \in \mathcal{E}_I : \bm{x}_J > \bm{0}, \bm{x}_{\setminus J} = \bm{0}\}. \end{align*}
   Note that $\mathcal{E}_I\neq \{\bm{x}_I : \bm{x} \in \mathcal{E}\}$. However, $\mathcal{E}^J_I=\{\bm{x}_I : \bm{x} \in \mathcal{E}^J\}$ for $\emptyset \neq J \subset I \subset \{1,\dots,d\}$ and $\mathcal{E}=\mathcal{E}_{\{1,\dots,d\}}$ and $\mathcal{E}^J=\mathcal{E}^J_{\{1,\dots,d\}}$ for $\emptyset \neq J \subset \{1,\dots,d\}$.
%%  and $\Lambda_I(S)=\Lambda(S \times [0,\infty)^{d-|I|})$.

    \paragraph{Lemma/Definition.} \textit{The following are equivalent and correspond to the traditional definition of extremal independence of $\bm{X}_A$ and $\bm{X}_C$.
    \begin{enumerate}[(i)]
      \itemsep0mm
    \item The exponent measure $\Lambda $ places mass only on $\mathcal{E}^I$ for which $I \subset A$ or $I \subset C$.
    \item $\Lambda(\bm{x}) = \Lambda_A(\bm{x}_A) + \Lambda_C(\bm{x}_C)$ for all $\bm{x} > 0$.
     \item $\Lambda_I(\mathcal{E}^I_I)=0$ if $I \cap A \neq \emptyset$ and $I \cap C \neq \emptyset$.
    \end{enumerate}
    }

  \begin{proof}[Proof of Lemma] For  $\bm{x} > 0$ all expressions in (ii) are finite and so (ii) is in fact equivalent to $\Lambda(\mathcal{S}^{A,C}(\bm{x})) = 0$ for $\bm{x} > 0$, where 
    \begin{align*}
      \mathcal{S}^{A,C}(\bm{x}) = \mathcal{E}_A \setminus [\bm{0},\bm{x}_A] \times \mathcal{E}_C \setminus [\bm{0},\bm{x}_C] = \{ \bm{z} \in \mathcal{E} : \exists{}\, a \in A: z_a > x_a, \exists{} \,c \in C: z_c > x_c\}.
      \end{align*}
    It is easy to check that $\mathcal{S}^{A,C}({\bm{x}}) \cap \mathcal{E}^I=\emptyset$ for $I \subset A$ or $I \subset C$.  Hence (i) implies (ii). Further, let $I \cap A \neq \emptyset$ and $I \cap C \neq \emptyset$ and note that $\mathcal{E}^I_I = \bigcup_{n \in \NN} \mathcal{E}^I_I(n)$ for $\mathcal{E}^I_I(n)=\{ \bm{z} \in \mathcal{E}^I_I: \bm{z} > \bm{n}^{-1}\}$. It is readily checked that $\mathcal{E}^I_I(n) \times [{0},{\infty})^{d-|I|} \subset \mathcal{S}^{A,C}(\bm{n}^{-1})$. Hence (ii) implies $\Lambda_I(\mathcal{E}^I_I(n))=\Lambda(\mathcal{E}^I_I(n) \times [{0},{\infty})^{d-|I|})=0$, and thus $\Lambda_I(\mathcal{E}^I_I)=\lim_{n \to \infty} \Lambda_I(\mathcal{E}^I_I(n)) = 0$, that is (iii). Finally (iii) implies (i), since  $\Lambda(\mathcal{E}^I) = \Lambda(\mathcal{E}^I_I \times \{\bm{0}\})\leq \Lambda(\mathcal{E}^I_I \times [{0},{\infty})^{d-|I|})=\Lambda_I(\mathcal{E}^I_I)$. \end{proof}

  \begin{remark}
    By the continuity of the (max-stable) probability measures with distribution functions $\exp(-\Lambda_I(\bm{x}))$ statement (ii) is also equivalent to
    \textit{
    \begin{enumerate}
      \item[(ii)'] $\Lambda(\bm{x}) = \Lambda_A(\bm{x}_A) + \Lambda_C(\bm{x}_C)$ for all $\bm{x} \geq \bm{0}$.
    \end{enumerate}
    }
  \end{remark}

  \enlargethispage{1cm}
  \begin{proof}[Proof of Theorem]

    The simple direction is ``old $\Rightarrow$ new'':
    Note that $\bigcup_{\emptyset \neq I \subset J}\mathcal{E}^I=\{\bm{x} \in \mathcal{E} : \bm{x}_{\setminus J} = 0\}$ for $J \in \{A,C\}$. Let $k \in A$. Recall that the law of $\bm{Y}^k$ equals $\Lambda$ restricted to $\mathcal{L}^k$.
    Since $\mathcal{L}^k \cap \bigcup_{\emptyset \neq I \subset C}\mathcal{E}^I = \emptyset$, this measure places mass only on $\mathcal{L}^k \cap \bigcup_{\emptyset \neq I \subset A}\mathcal{E}^I = \{\bm{x} \in \mathcal{L}^k, \bm{x}_C=0\}=\mathcal{L}^k_A \times \{\bm{0}\}$ and therefore splits as the product measure ``$\Lambda_A$ restricted to $\mathcal{L}^k_A$ times a point mass at $\bm{0}$'', as desired. The same holds true with reversed roles for $A$ and $C$.
    
    The converse ``new $\Rightarrow$ old'' requires slightly more work: If $k \in I$ we know that the law of $\bm{Y}^k_I$ is $\Lambda_I$ restricted to $\mathcal{L}^k_I$. Hence, $\bm{Y}_A \perp_e \bm{Y}_C$  implies that the law of $\bm{Y}^k=(\bm{Y}_A,\bm{Y}_C)|Y_k > 1$ (which is $\Lambda$ restricted to $\mathcal{L}^k$) splits into a product measure with first component $\Lambda_A$ on $\mathcal{L}_A^k$ if $k \in A$ and second component $\Lambda_C$ on $\mathcal{L}_C^k$ if $k \in C$. For $\bm{x} > \bm{1}$, the set $[\bm{x},\bm{\infty})$ lies in the intersection of all $\mathcal{L}^k$, $k \in A \cup C$. Hence $\Lambda([\bm{x},\bm{\infty}))=\Lambda_A([\bm{x}_A,\bm{\infty}))\mu_C([\bm{x}_C,\bm{\infty}))=\mu_A([\bm{x}_A,\bm{\infty}))\Lambda_C([\bm{x}_C,\bm{\infty}))$ for some measures $\mu_A$ and $\mu_C$. Suppose there is some $\bm{x}^* > \bm{1}$, for which $\Lambda([\bm{x^*},\bm{\infty})) \neq 0$. Then it follows that $\Lambda([\bm{x},\bm{\infty}))$ is a constant multiple of $\Lambda_A([\bm{x}_A,\bm{\infty}))\Lambda_C([\bm{x}_C,\bm{\infty}))$ for $\bm{x} > \bm{1}$. However, this contradicts the homogeneity of $\Lambda,\Lambda_A,\Lambda_C$, which all satisfy $\mu(tS)=t^{-1}\mu(S)$ for $S$ bounded away from $\bm{0}$ (in their respective spaces). Hence $\Lambda([\bm{x},\bm{\infty}))=0$ for all $\bm{x} > \bm{1}$ and, by the homogeneity of $\Lambda$, for all $\bm{x} > \bm{0}$. This shows that the exponent measure $\Lambda$ cannot place mass on $\mathcal{E}^{A \cup C}=\mathcal{E}^{\{1,\dots,d\}}$. By moving to lower dimensional margins, the same argument can be repeated to show that $\Lambda_I(\mathcal{E}^I_I)=0$, whenever $I$ contains elements from both $A$ and $C$. 
    \end{proof}

\end{document}